\documentclass{amsart}

\usepackage{amsmath}
\usepackage{amsfonts,amssymb}
\usepackage{graphics}
\begin{document}

\newtheorem{theorem}{Theorem}
\newtheorem{corollary}[theorem]{Corollary}
\newtheorem{thm}{Theorem}[section]
\newtheorem{cor}[thm]{Corollary}
\newtheorem{lem}[thm]{Lemma}
\newtheorem{prop}[thm]{Proposition}
\theoremstyle{definition}
\newtheorem{definition}[thm]{Definition}
\theoremstyle{remark}
\newtheorem{rem}[thm]{Remark}
\newtheorem{examp}{Example}
\numberwithin{equation}{section}

\newcommand{\R}{\mathbb{R}}
\newcommand{\ZM}{\mathbb{Z}}
\newcommand{\QM}{\mathbb{Q}}
\newcommand{\NM}{\mathbb{N}}
\newcommand{\CM}{\mathbb{C}}

\newcommand{\eps}{\varepsilon}

\newcommand{\Z}{\mathbb{Z}}
\newcommand{\conj}[1]{\overline#1}
\newcommand{\diag}{\rm diag}
\newcommand{\disps}{\displaystyle}
\newcommand{\ep}{\qed\endtrivlist}
\newcommand{\eqskip}{ \vspace*{2mm} }
\newcommand{\inta}{ {\disps \int_{A}} }
\newcommand{\ko}{K_{0}}
\newcommand{\mfam}{ \mathcal{M} }
\newcommand{\nsu}{ \mathcal{N}_{u} }
\newcommand{\om}{ \Omega }
\newcommand{\parag}{ \hspace{-0,65cm} }
\newcommand{\sca}{ \mathcal{S} }
\newcommand{\ub}{\bar{u}}
\newcommand{\vphi}{\varphi}

\title[Nodal lines of the Laplacian on surfaces]{Closed nodal lines
and interior hot spots of the second eigenfunction of
the Laplacian on surfaces}
\author{Pedro Freitas}
\thanks{Partially supported by FCT, Portugal}
%
%
\address{Departamento de Matem\'{a}tica,
Instituto Superior T\'{e}cnico, Av.Rovisco Pais, 1049-001 Lisboa,
Portugal.} \email{pfreitas@math.ist.utl.pt}
\date{\today}
\begin{abstract}
We build a one--parameter family of $S^{1}-$invariant metrics on the
unit disc with fixed total area for which the second eigenvalue of the
Laplace operator in the case of both Neumann and Dirichlet boundary
conditions is simple and has an eigenfunction with a closed nodal
line. In the case of Neumann boundary conditions, we also prove that
this eigenfunction attains its maximum at an interior point, and thus
provide a counterexample to the hot spots conjecture on a simply
connected surface. This is a consequence of the stronger result that
within this family of metrics any given (finite) number of
$S^{1}-$invariant eigenvalues can be made to be arbitrarily small,
while the non--invariant spectrum becomes arbitrarily large.
\end{abstract}

\maketitle
\section{Introduction}

A conjecture of J. Rauch from 1974 states that the eigenfunction
corresponding to the second eigenvalue of the Laplace operator on a
domain with Neumann boundary conditions attains its maximum and
minimum on the boundary. As has been pointed out in~\cite{babu}, for
instance, this is related to the location of the points of maxima of
solutions of the heat equation (hot spots) and it is basically
equivalent to saying that, for {\it most} initial conditions, these
hot spots move towards the boundary as time goes to infinity.

Recently, Burdzy and Werner~\cite{buwe} gave a counterexample to this
conjecture on a domain in $\R^{2}$ with two holes and posed the
question of whether it would still be possible to find a
counterexample on a doubly connected domain, or whether the conjecture
would hold in that case -- note that it is not known if the conjecture
holds on simply connected domains.

On the positive side, Kawohl~\cite{kawo} has shown that for the case
of domains of the form $D\times (a,b)$, where $D\subset\R^{n-1}$ has a
$\mathcal{C}^{0,1}$ boundary, the values attained by a second
eigenfunction on the domain are less than or equal to the values it
takes on the boundary. More recently, Ba\~{n}uelos and Burzdy proved
that the conjecture holds in the case of some special domains which
include, among others, convex domains with a line of
symmetry~\cite{babu}.

This problem is related (although not necessarily
equivalent) to the nonexistence of closed nodal lines of the second
eigenfunction, where the nodal set of an eigenfunction is defined to
be the closure of the subset of the domain where the eigenfunction vanishes.
With the exception of Courant's nodal domain theorem which states that the
nodal set of a $k^{\rm th}$ eigenfunction divides the domain into at
most $k$ subregions~\cite{cohi}, very little is known about the
general structure of such sets -- see~\cite{dofe,grje,jeri2}.
Courant's result implies that a second eigenfunction always divides
the domain into exactly two subregions and a conjecture of Payne from
1967 for the case of Dirichlet boundary conditions states that any
such eigenfunction cannot have a closed nodal line~\cite{payn1}. This
was proved by Payne in 1973 under some symmetry and convexity
assumptions on the domain~\cite{payn2}. Within the last ten years
there have been several developments regarding the existence or not of
closed nodal lines of the second eigenfunction for the Dirichlet
Laplacian -- see~\cite{ales,four,hhn,jeri,jeri2,lin,lini,mela,putt}.
Although it has been shown that the conjecture holds for convex
domains~\cite{ales,mela}, again it is not known whether or not there
exist simply connected sets where the second eigenfunction will have a
closed nodal line. So far the only known counterexample requires the
boundary of the domain to have at least three components~\cite{hhn}.
Note that, in the case of the Neumann problem, it is not too difficult 
to prove that on simply connected domains the second eigenfunction 
cannot have a closed nodal line~\cite{band,kawo,payn1}.

The question of whether or not results similar to these hold in the
case of manifolds with boundary has been raised by S.T. Yau -- see,
for instance, Problem 45 in the Chapter {\it Open problems in
differential geometry} in~\cite{scya}. Note that in this case it is
not difficult to find a counterexample to the closed nodal line
conjecture if one considers doubly connected surfaces. For this, it is
sufficient to think of long cylinders, for which the nodal line is a
circle located halfway between the top and the bottom. Incidentally,
this also justifies the statement that the hot spots and the closed
nodal line problems are not necessarily equivalent, as in this case
the maximum and minimum are still attained at the boundary. As far as
we know, the only results for the case of surfaces are given in the
book by Bandle~\cite{band}, where it is shown that in the Neumann case
and under some conditions on the metric the second eigenvalue is not
simple and the maximum and minimum are attained at the boundary. It
is stated in that book that one cannot expect for this to
be true in general, but no counterexample is given.

The main purpose of this note is to present examples of
simply--connected surfaces with boundary for which both the closed
nodal line and the hot spots conjectures fail. In other words, for
these surfaces the eigenfunctions associated with the second
eigenvalue of the Laplace operator in both the Neumann and the
Dirichlet case have closed nodal lines which do not touch the
boundary. Furthermore, in the Neumann case, this eigenfunction can be
chosen in such a way that it attains its maximum at an interior point.
This shows that simple connectivity by itself is not a sufficient 
condition for the above conjectures to hold.

These counterexamples are corollaries to the stronger result which
asserts the existence of a family of $S^{1}$-invariant metrics with
constant area on the unit disc for which the first $m$ eigenvalues are
simple, for any given positive integer $m$. For $S^{1}$-invariant metrics,
the spectrum can be divided into an invariant part corresponding to
the spectrum of the Laplacian acting on $S^{1}$-invariant functions,
and a non--invariant part. The counterexamples given here are then
made possible by the fact that we are able to choose a family of
metrics for which the invariant and the non--invariant parts of the
spectrum may be {\it separated}. More precisely, for any positive
integer $m$ we have that the first $m$ invariant eigenvalues can be
made to be arbitrarily small, while the first non--invariant
eigenvalue becomes arbitrarily large.

\section{Main result}

The main result from which the counterexamples are then a
straightforward consequence is the following
\begin{theorem}\label{t1}
There exists a one--parameter family $\mfam$ of smooth $S^{1}-$in\-variant
metrics on the unit disc with positive curvature and total area $\pi$, such
that given any positive integer $m$ and real number $\eps$ there exists a
subset $\mfam_{\eps}$ of $\mfam$ with $\eps$ on an open interval for which
the first $m$ $S^{1}-$invariant
eigenvalues of the Laplace operator both with Neumann and Dirichlet
boundary conditions are smaller than $\eps$. On the other hand, the
non--invariant spectrum remains uniformly bounded away from zero in
$\mfam$ and becomes arbitrarily large in $\mfam_{\eps}$ as $\eps$ goes
to zero.
\end{theorem}
Since the Gaussian curvature is positive, these surfaces can actually 
be isometrically embedded in $\R^{3}$~\cite{nire}. As we shall see
in Section~\ref{prt1}, the curvature of such metrics can be made to be
arbitrarily close to zero except at the centre of the disc.

From the proof of the theorem, it follows that the first $m$ invariant
eigenvalues must be simple, and, in the case of the second Neumann
eigenvalue, that the corresponding (invariant) eigenfunction is
strictly monotone along radial lines and changes sign. We thus obtain
the following
\begin{corollary}
Given any positive integer $m$ there exists a family $\mfam$ of
$S^{1}-$invariant metrics on the unit disc with positive curvature and
total area $\pi$, for which the first $m$ eigenvalues of the Laplace
operator with both Neumann and Dirichlet boundary conditions are
simple. In both cases, the eigenfunctions corresponding to the $j^{\rm
th}$ eigenvalue ($j=2,\ldots,m$) are also invariant by $S^{1}$ and
have $j-1$ nodal lines which are closed disjoint circumferences
dividing the disc into $j$ nodal domains. In the Neumann case the
second eigenfunction may be chosen in such a way that its (strict)
maximum is attained at the origin.
\end{corollary}
The estimates obtained in the proof allow us to write a more
quantitative version of this result which we mention in
Section~\ref{prt1}.

Counterexamples to this type of conjectures have usually been obtained
by exhibiting a domain for which the second eigenvalue can be proven
to be simple~\cite{buwe,hhn,lini}. If the domain has some symmetry, it
will then be inherited by the eigenfunction and hence also by the nodal set.
In order to prove the above results, we shall use a variation of this
technique and consider one--parameter families of metrics on the disc
for which we are able to prove that the $m^{\rm th}$ invariant
eigenvalue is smaller than the first non--invariant eigenvalue.

A fundamental (standard) ingredient in the proof is the reduction of
the original two--dimensional problem to a sequence of
one--dimensional problems, which we achieve by means of standard polar
coordinates. Since we are interested in estimates for higher
eigenvalues, we then perform a change of variables in order to avoid
function weights in the orthogonality conditions which appear in the
corresponding variational formulations
-- see Section~\ref{abst}. At this point we should remark that it
is also possible to proceed by means of a different technique using
symplectic coordinates. An example of this can be found
in~\cite{abfr}, where the behaviour of the invariant spectrum
for $S^{1}-$invariant metrics on $S^{2}$ was studied. We shall briefly
indicate in Section~\ref{abst} how these two different coordinate systems
are related in this case.

Regarding the choice of metrics, we point out that with the proper
parame\-trization in isothermic coordinates (see Section~\ref{abst}),
the eigenvalue problem on the surface becomes equivalent to that on a
flat disc of inhomogeneous density, and so our results also apply in
that case. On the other hand, this suggests that the we build the
family $\mfam$ by making the density much higher close to the centre
than near the boundary, so that there is a strong resistance to the
movement of the hot spots towards the outside regions of the disc. We
remark that in the opposite case where the (radially symmetric)
density increases as we move away from the centre, it is known that
the hot spots conjecture holds~\cite{band}.

\section{Abstract surfaces\label{abst}}

In this section we collect the main facts about abstract surfaces that
will be needed in the paper. We shall follow the exposition in the
book by Bandle~\cite{band} very closely, and begin by introducing the
concept of an abstract surface $\sca$. We then derive the expressions
that will be used in the sequel, which include the Laplace--Beltrami
operator in $\sca$ given in conformal coordinates and the variational
formulation for the problem. We shall also indicate the expression for
the Gaussian curvature of $\sca$.


Let $D$ be a domain in the $(x,y)-$para\-meter plane, and
$d\sigma^{2}$ the Riemannian metric in $D$ defined by the quadratic form
\[
d\sigma^{2}=E(x,y)dx^{2}+2F(x,y)dxdy+G(x,y)dy^{2}.
\]
\begin{definition}
A domain $D\subset\R^{2}$ with the Riemannian metric $d\sigma$ is
called an abstract surface and will be denoted by $\sca=(D,d\sigma)$.
$\sca$ is said to be in its isothermic (or conformal) representation
if $d\sigma^{2}=p(x,y)ds^{2}$ ($E=G=p$ and $F=0$), where $ds$ denotes
the linear element of the Euclidean plane.
\end{definition}

The eigenvalue problems that we are interested in are thus
\begin{equation}
\label{lapsurf}
\begin{array}{rl}
\Delta_{\sca}u+\gamma u = 0 & \mbox{in } D,\eqskip\\
\frac{\disps \partial u}{\disps \partial \nu} = 0 \mbox{ or } u=0 &
\mbox{on } \partial D,
\end{array}
\end{equation}
where $\Delta_{\sca}$ denotes the Laplace--Beltrami operator on $\sca$
and $\nu$ the conormal derivative. We shall denote the eigenvalues of
the Neumann and Dirichlet problems by $\mu_{j}$ and $\lambda_{j}$,
$j=1,\ldots$, respectively, and always assume that they are written in
increasing order.

The Laplace--Beltrami operator is now given by
\[
\Delta_{\sca} = \frac{\disps 1}{\disps W}\left[
\frac{\disps\partial}{\partial x}\left(\frac{ \disps G
\frac{\disps\partial}{\disps\partial x}-F
\frac{\disps\partial}{\disps\partial y}}{\disps W}\right) +
\frac{\disps\partial}{\partial y}\left(\frac{ \disps E
\frac{\disps\partial}{\disps\partial y}-F
\frac{\disps\partial}{\disps\partial x}}{\disps W}\right)
\right],
\]
where $W=\sqrt{EG-F^{2}}$. In the case where $\sca$ is in its
isothermic representation, this expression simplifies to
\[
\Delta_{\sca} = \frac{\disps 1}{\disps p(x,y)}\Delta,
\]
where $\Delta$ now denotes the usual Laplacian operator in $\R^{2}$.
This means that the eigenvalue problem~(\ref{lapsurf}) on $\sca$
becomes
\begin{equation}
\label{lapiso}
\begin{array}{rl}
\frac{\disps 1}{\disps p(x,y)}\Delta u+\gamma u = 0 & \mbox{in } D,\eqskip\\
\frac{\disps \partial u}{\disps \partial \nu} = 0 \mbox{ or } u=0 &
\mbox{on } \partial D,
\end{array}
\end{equation}
As mentioned in the Introduction, this is equivalent to the problem of
an inhomogeneous membrane whose density is given by the function $p$.

Finally, we note that in this case, the area of $\sca$ and the
Gaussian curvature are given by
\[
A_{p}(\sca) = \int_{D}p(x,y)dxdy\
\mbox{ and }\
K = -\frac{\disps 1}{\disps 2p} \Delta\left[\log(p)\right],
\]
respectively. The expression for the total curvature is then
\[
\omega(\sca) = -\frac{\disps 1}{\disps 2}{\disps \int_{D}} \Delta\left[\log(p(x,y))\right]dxdy.
\]


From this point on we shall concentrate on the Dirichlet problem, and
mention the necessary changes for the Neumann case in
Section~\ref{prt1}.

We consider the eigenvalue problem~(\ref{lapiso}) in the case where
the function $p$ is radially symmetric, that is, $p=p(r)$. In polar
coordinates and after separation of variables, the eigenvalue
problem~(\ref{lapiso}) then reduces to the sequence of
one--dimensional eigenvalue problems
\begin{equation}
\label{seq1d}
(r \vphi')' + \left[\lambda rp(r) - \frac{\disps k^{2}}{\disps
r}\right]\vphi
= 0, \ \ r\in(0,1),
\end{equation}
where $k^{2}$, $k=0,1,\ldots$, are the eigenvalues of the Laplacian on
the circle. For $k=0$ we have Neumann boundary conditions at $0$ and
Dirichlet at $1$, that is $\vphi'(0)=\vphi(1)=0$. This corresponds to
the invariant spectrum of the original problem, and the associated
eigenfunctions are the radially symmetric functions $\vphi$. Note that
the $i^{\rm th}$ eigenfunction of the one--dimensional problem has
$i-1$ zeros on $(0,1)$, and thus the corresponding eigenfunction of
the original problem divides the disc into $i$ nodal domains.

For positive values of $k$ we have Dirichlet boundary conditions at
both ends of the interval, and this now gives the non--invariant part
of the spectrum which is made up of eigenvalues with multiplicity two.
In this case, two linearly independent eigenfunctions are given by
$\vphi(r)\cos(k\theta)$ and $\vphi(r)\sin(k\theta)$.

To obtain the variational formulation corresponding to the eigenvalue
problems above we consider the space $\mathcal{C}^{1}(0,1)$ with the
inner product defined by
\[
(\vphi,\psi) = \int_{0}^{1}rp(r)\vphi(r)\psi(r) + r\vphi'(r)\psi'(r) dr.
\]
Let now $H^{1}_{0}$ be the Sobolev space which is obtained as the
closure of $\mathcal{C}^{\infty}_{0}(0,1)$ with respect to the norm
induced by the above inner product. Then the eigenvalues of the
spectral problem~(\ref{seq1d}) when $k=0$ are given by
\begin{equation}
\label{varfor}
\lambda^{0}_{j}=\lambda^{0}_{j}(p) = \inf_{\vphi\in H^{1}_{0,j}}\frac{\disps
\int_{0}^{1}r[\vphi'(r)]^{2}dr}
{\disps \int_{0}^{1}rp(r)\vphi^{2}(r)dr}, \ \ j=1,\ldots
\end{equation}
where $H^{1}_{0,1}$ is the closure with respect to the above norm of
the space of $\mathcal{C}^{\infty}$ functions on $[0,1)$ with compact
support on this interval,
\[
H^{1}_{0,j} = H^{1}_{0,j-1} \cap
\left\{\vphi: \int_{0}^{1}rp(r)\vphi_{j-1}(r)\vphi(r)dr = 0\right\},
\ \ j=2,\ldots,
\]
and $\vphi_{j}$ is the eigenfunction corresponding to
$\lambda^{0}_{j}$.

In the case of the eigenvalues for each of the remaining problems
we have
\begin{equation}
\label{varforb}
\lambda^{k}_{j}=\lambda^{k}_{j}(p) = \inf_{\vphi\in H^{1}_{0,j}}
\frac{\disps \int_{0}^{1}r[\vphi'(r)]^{2}dr +
k^{2}\int_{0}^{1}\frac{\disps \vphi^{2}(r)}{\disps r}dr} {\disps
\int_{0}^{1}rp(r)\vphi^{2}(r)dr}, \ \ j=1,\ldots, \  k=1,\ldots,
\end{equation}
and where now the spaces $H^{1}_{0,j}$ are defined in a similar way 
as above but starting with $H^{1}_{0}$.

As we are interested in higher eigenvalues, we have to take into
account the orthogonality conditions which appear in the definition of
the spaces $H^{1}_{0,j}$. Since these conditions become much simpler
to handle if the inner product considered has a constant weight
function instead of $rp(r,\delta)$, in order to show that
$\lambda^{0}_{m}(p)$ is arbitrarily small within a family of metrics
it is convenient to make a change of variables in the variational
formulation of the problem. Thus, we want to find a function $r=r(z)$
and a constant $c$ such that
\[
r(z)p\left(r(z)\right)r'(z) = c
\]
and the interval $(0,1)$ is mapped onto $(0,1)$. This will be the case
if we take, for instance,
\begin{equation}
\label{change}
z(r) = \frac{\disps 1}{\disps c}{\disps \int_{0}^{r}}sp(s)ds \ \
\mbox{and} \ \ c = {\disps \int_{0}^{1}sp(s)ds}=\frac{\disps
A_{p}(\sca)}{\disps 2\pi}.
\end{equation}
The Raleigh quotients in~(\ref{varfor}) and~(\ref{varforb}) then become
\begin{equation}
\label{varforc}
\frac{\disps 4\pi^{2}}{\disps A_{p}^{2}(\sca)}
\frac{\disps
\int_{0}^{1}r^{2}(z)p\left(r(z)\right)[\psi'(z)]^{2}dz}
{\disps
\int_{0}^{1}\psi^{2}(z)dz}+
k^{2}\frac{\disps \int_{0}^{1}\left[r^{2}(z)p\left(r(z)\right)
\right]^{-1}[\psi(z)]^{2}dz}
{\disps\int_{0}^{1}\psi^{2}(z)dz},
\end{equation}
for $k=0,\ldots$. This transformation fixes both endpoints of the interval,
and the inner product is changed accordingly to
\[
(\phi,\psi) = A_{p}(\sca){\disps \int_{0}^{1}} \phi\psi dz +
\frac{\disps 4\pi^{2}}{\disps A_{p}(\sca)}{\disps \int_{0}^{1}}
r^{2}(z)p(r(z))\phi'\psi'dz.
\]
We thus have that the orthogonality conditions used to define the spaces
appearing in the variational problem are now formed with respect to the usual
inner product. We shall use the same notation for these spaces as above, as it
will be clear from the context which inner product is being used.

We remark that it is also possible to choose the change of variables
above in such a way that the problem becomes equivalent to that
obtained in~\cite{abfr} using symplectic coordinates. We thus obtain a
relation between symplectic coordinates and isothermic coordinates for
this particular case. To do this, we need the relation between
$r$ and $z$ to be such that
\begin{equation}
\label{pimplicit}
r^{2}p(r) = \frac{\disps 1}{\disps g(z)} \mbox{ and } r'(z) = rg(z).
\end{equation}
Here $g:(-1,1)\to(0,+\infty)$ is a function of the form
\[
g(z) = \frac{\disps 1}{1-z^{2}} + h(z),
\]
where $h$ is a $\mathcal{C}^{\infty}[-1,1]$ function such that $g>0$,
and it has been assumed that the area of the full surface is $4\pi$.
We then obtain that the change of variables defined by
\[
r(z) = \sqrt{\frac{\disps \eps}{\disps 2-\eps}} e^{-\disps\int_{z}^{1}
h(t)dt}\sqrt{\frac{\disps 1+z}{\disps 1-z}}\
\ \ \eps>0,
\]
transforms~(\ref{varfor}) into
\[
\frac{\disps \int_{-1}^{1-\eps}\frac{\disps 1}{\disps g(z)}[\psi'(z)]^{2}dz}
{\disps \int_{-1}^{1-\eps}\psi^{2}(z)dz}.
\]
The function $p$ is then given from~(\ref{pimplicit}) by
\[
p(r) = \frac{\disps 1}{\disps r^{2}g(z(r))},
\]
although it should be pointed out that in general it will not be
possible to obtain an explicit expression for it in terms of $r$. Note
that because polar coordinates have a singularity at the origin, while
symplectic coordinates have two singularities at both the North and South 
poles, we have to restrict ourselves to the interval
$(-1,1-\eps)$. In this case, the change of variables given 
by~(\ref{change}) has the advantage that it keeps the interval fixed.

As an example, consider the case of the standard sphere for which we
have that $h$ is the zero function and then
\[
p(r) = \frac{\disps 4\eps(2-\eps)}
{\disps\left[\eps+(2-\eps)r^{2}\right]^{2}}
\]
as expected.

\section{Proof of Theorem 1\label{prt1}}

The first $m$ eigenvalues of the original problem will be simple if
and only if we have that $\lambda^{0}_{m}$ is smaller than
$\lambda^{k}_{1}$ for all positive $k$. This follows from the fact
that the invariant spectrum is the spectrum of a second order ordinary
differential operator of the type in~(\ref{seq1d}) with $k=0$ and thus
all its eigenvalues are simple. Clearly it is enough to ensure the
above condition for $k=1$, since $\lambda^{k}_{1}<\lambda^{k+1}_{1}$
for all $k$. To do this, we consider the one--parameter family of
smooth metrics corresponding to
\[
p(r,\delta) = \frac{\disps \alpha}{\disps r^{2}+\delta}, \ \ \
\ \alpha=
\frac{\disps 1}{\disps \log\left(\frac{\disps 1+\delta}{\disps \delta}\right)},
\]
for positive $\delta$. The parameter $\alpha$ is included as a
normalizing factor so that all elements in this family of metrics
have fixed area equal to $\pi$.

Using the expression given in Section~\ref{abst} we easily obtain
\[
K(r,\delta) = \frac{\disps 2\delta}{\disps r^{2}+\delta}
\log\left(\frac{\disps 1+\delta}{\disps \delta}\right),
\]
and thus the curvature is positive and converges to $0$ on $(0,1]$ as
$\delta$ goes to $0$. On the other hand, the total curvature is given
by
\[
\omega(\delta) = \frac{\disps 2\pi}{\disps 1+\delta}.
\]


Using the change of variables given by~(\ref{change}) we obtain
\[
r = \sqrt{\delta\left(e^{z/\alpha}-1\right)},
\]
and the variational problem~(\ref{varforc}) becomes
\begin{equation}
\label{varform2}
\begin{array}{l}
\lambda^{k}_{j}=\lambda^{k}_{j}(p) = {\disps\inf_{\psi\in
H^{1}_{0,j}}}\left(4\alpha\frac{\disps
\int_{0}^{1}\left(1-e^{-z/\alpha}\right)[\psi'(z)]^{2}dz}
{\disps \int_{0}^{1}\psi^{2}(z)dz}+\right.\eqskip \\
\hspace*{2cm}\left.+ \frac{\disps k^{2}}{\disps \alpha} \frac{\disps
\int_{0}^{1}\frac{\disps 1}{\disps \left(1-e^{-z/\alpha}\right)}\psi^{2}(z)dz}
{\disps \int_{0}^{1}\psi^{2}(z)dz}\right), \ \ j=1,\ldots, \ \ k=0,\ldots.
\end{array}
\end{equation}

We begin by noting that since $\alpha$ is positive the second term is 
always greater than or equal to $k^{2}/\alpha$ so that we have
\[
\lambda_{1}^{1}(p) \geq \frac{\disps 1}{\disps \alpha},
\]
which becomes unbounded as $\delta$ goes to zero.

We shall now show that any finite number of invariant
eigenvalues can be made to be arbitrarily small. We have that
\[
\frac{\disps
\int_{0}^{1}\left(1-e^{-z/\alpha}\right)[\psi'(z)]^{2}dz}
{\disps \int_{0}^{1}\psi^{2}(z)dz} \leq
\frac{\disps
\int_{0}^{1}[\psi'(z)]^{2}dz}{\disps \int_{0}^{1}\psi^{2}(z)dz}
\]
and so by the monotonicity principle derived from Poicar\'{e}'s 
principle it follows that
\[
\lambda_{j}^{0}(p)\leq \alpha (2j-1)^{2}\pi^{2}, \ \ j=1,\ldots,
\]
and thus, given a fixed eigenvalue $\lambda^{0}_{j}(p)$,
\[
\lim_{\delta\to0^{+}} \lambda^{0}_{j}(p) = 0.
\]

The estimates obtained above allow us to conclude that if
\[
\delta<\frac{\disps 1}{\disps e^{(2j-1)\pi}-1},\ \ j=1,\ldots,
\]
then the $j^{\rm th}$ invariant eigenvalue is smaller than the first
non--invariant eigenvalue. This is a very rough estimate,
since the logarithmic term is taking into account the asymptotic
behaviour of the eigenvalues, while the crossing between
the two eigenvalues as the parameter $\delta$ decreases must occur at
a point larger than $1$. Similar results can be obtained for higher
eigenvalues.


Due to the fact that the Dirichlet and Neumann eigenvalues satisfy
$\mu_{j}\leq\lambda_{j}$ for all integer $j$, the proof for the
Neumann problem follows from the Dirichlet case. It is also possible
to obtain some independent estimates by following essentially along
the same lines as in the Dirichlet case, with only some minor changes.
These have to do with the spaces considered, which should now be
\[
\begin{array}{lcl}
H^{1}_{1} = H^{1} & \mbox{and} & H^{1}_{j} = H^{1}_{j-1} \cap \left\{
\vphi: \int_{0}^{1}\vphi_{j-1}(z)\vphi(z)dz = 0\right\}
\end{array}
\]
in the case of invariant eigenvalues, and $H^{1}_{0}$ in~(\ref{varforb}) is
now the closure of the space of $\mathcal{C}^{\infty}$ functions with compact
support on $(0,1]$, with respect to the same norm
as before. We may then proceed in a similar fashion to obtain the following
estimate for the $j^{\rm th}$ invariant Neumann eigenvalue:
\[
\mu^{0}_{j}(p)\leq 4\alpha(j-1)^{2}\pi^{2}.
\]

From~(\ref{seq1d}) with $k=0$ we see that any solution $\vphi$ of that
equation with precisely one zero on the interval $(0,1)$ and with
$\vphi(0)$ positive will be strictly decreasing on this interval.
Hence its maximum is attained at the origin.

\subsection*{Acknowledgments}
This work was carried out while I was visiting the Department of
Mathematics of the Royal Institute of Technology in Stockholm, Sweden.
I would like to thank the people there and, in particular, Ari Laptev,
for their hospitality. I would also like to thank Miguel Abreu for the
many helpful conversations.


\end{document}